\documentclass[a4paper]{article}

\usepackage{amsmath}
\usepackage{amsfonts}
\usepackage{amssymb}         
\usepackage{amsopn}
\usepackage{theorem}          
\usepackage{latexsym}
\usepackage{graphicx}        
\usepackage{mathrsfs}           

\theoremstyle{changebreak}                

\newenvironment{proof}
 {{\sl Proof.}\hspace*{1 ex}}%
 {{\nopagebreak\hspace*{\fill}$\Box$\par\vspace{12pt}}}
\begin{document}

\thispagestyle{empty}
\begin{center} 

{\LARGE On a Linear Property of Bilinearly Defined Sets in $\mathbb{R}^n$}
\par \bigskip
{\sc Leo Liberti\footnote{E-mail: {\tt l.liberti@ic.ac.uk}.}}
\par  
{\it Centre for Process Systems Engineering \\
Imperial College of Science, Technology and Medicine \\
London SW7 2BY \\
United Kingdom}
\par \medskip 9 July 2002
\end{center}
\par \bigskip

\begin{abstract} 
Certain types of bilinearly defined sets in $\mathbb{R}^n$ exhibit a higher 
degree of linearity than what is apparent by inspection. 
\end{abstract}

\ \hfill \hrule \hfill \ 

Let $n,m\in\mathbb{N}$ such that $n<m$; $w,x\in\mathbb{R}^n$;
$y\in\mathbb{R}$; $J$ be an index set of size $n-m$ such that for
all $j\in J$ we have $j\le n$; $A=(a_{ij})$ an $m\times n$ matrix
having rank $m$; $b\in\mathbb{R}^m$;
\begin{eqnarray}
  \label{C}
  C&=&\{ (x,w,y) \; | \; Ax = b \wedge \forall j\le n \ (w_j=x_jy) \}; \\
  R_J&=&\{ (x,w,y) \; | \; Ax = b \wedge Aw-by=0 \wedge \forall j\in J
       \ (w_j=x_jy) \}.
  \label{R}
\end{eqnarray}
Note that since $y\in\mathbb{R}$, the product $by$ is simply
$(b_1y,\ldots,b_my)^T$. We shall show that there is at least a $J$
such that $C=R_J$. This in effect means that $m$ of the $n$ bilinear
terms that define the set $C$ can be replaced by the linear system
$Aw-by=0$. Thus the set $C$ has a higher degree of linearity than what
is apparent by inspection of eqn. (\ref{C}).

\noindent
{\bf Theorem} \\
  $\exists J \ (C = R_J)$.

\begin{proof}
If $(x,w,y)\in C$, then by definition it satisfies $Ax=b$; multiplying
this system by $y$ we obtain $yAx-by=0$, and since $w_j=x_jy$ for all
$j\le n$ we have $Aw-by=0$. Hence $(x,w,y)\in R_J$ for all $J$. Now
because $\mbox{rk}(A)=m$, by applying Gaussian elimination we can find
a permutation $\pi\in S_n$ (the symmetric group of order $n$), a
vector $b^\prime\in \mathbb{R}^m$ and a matrix $A^\prime = (A^\prime_0
| A^\prime_1 )$, where:
\begin{eqnarray*}
A^\prime_0 = \left( \begin{array}{cccc} 
 a^\prime_{11}&a^\prime_{12}&\ldots&a^\prime_{1m}\\
 0            &a^\prime_{22}&\ldots&a^\prime_{2m}\\
 \vdots       &\ddots    &\ddots   &\vdots       \\
 0            &\ldots    &0     &a^\prime_{mm} \end{array} \right)
&,&
A^\prime_1 = \left( \begin{array}{ccc} 
 a^\prime_{1,m+1} & \ldots & a^\prime_{1n}\\
 a^\prime_{2,m+1} & \ldots & a^\prime_{2n}\\
 \vdots & \ddots & \vdots       \\
 a^\prime_{m,m+1} & \ldots & a^\prime_{mn} \end{array} \right),
\end{eqnarray*}
$A^\prime_0$ being an $m\times m$ square matrix of rank $m$, such that
$A^\prime \pi(x) = b^\prime$. Let $J^\prime=\{m+1,\ldots,n\}$,
$J=\pi^{-1}(J^\prime)$, and $(x,w,y)\in R_J$. Since $b=Ax$ we can
substitute it into the system $Aw-by=0$ to get $Aw-yAx=0$,
i.e. $A(w-xy)=0$. This is an underdetermined linear system with $m$
equations and $n$ variables
\begin{equation*}
  \forall i\le m \ \left( \sum_{j=1}^n a_{ij} (w_j-x_jy) = 0 \right),
\end{equation*}
which (by Gaussian elimination) has the same solution space as the
system $A^\prime\pi(w-xy)=0$. Since for all $j\in J$ we have, by
definition of $R_J$, $w_j=x_jy$, this is in turn equivalent to
\begin{equation}
  A^\prime_0 \bar{\pi}(w-xy) = 0 
  \label{pisystem}
\end{equation}
where $\bar{\pi}$ is a restriction of $\pi$ to $\{1,\ldots,
n\}\backslash J$ (so that $\bar{\pi}(j)=\pi(j)$ only if $j\not\in
J$). The system (\ref{pisystem}) is square with full rank $m$ and thus
has a unique solution $\bar{\pi}(w-xy) = 0$, i.e.
\begin{equation*}
  \forall j\not\in J \; (w_j=x_jy).
  \label{solution}
\end{equation*}
Hence $(x,w,y)$ satisfies the definition of $C$, i.e. $R_J\subseteq C$.
\end{proof}

The geometrical implications of this theorem are that the
intersection in $\mathbb{R}^n$ of a set of bilinear terms like those
described above is a hypersurface containing a degree of
linearity. By exploiting this linearity we are able to replace some of
the bilinear terms with linear terms. This fact is useful in a number
of applications, e.g. linearly constrained bilinear programming.


\bibliographystyle{alpha}
\bibliography{phd}

\end{document}